\documentclass[twoside,a4paper]{amsart}

\PassOptionsToPackage{pdfauthor={Vladimir V. Kisil},%
    pdftitle={Cycles cross ratio: an invitation},%
    pdfsubject={mathematics},%
    backref=page,%
  pdfkeywords={symmetries}
}{hyperref}
\IfFileExists{eulervm.sty}{%

  \usepackage{hyperref}
\usepackage{graphicx}
\usepackage{eulervm}
}
{
}
\providecommand{\cites}[1]{\cite{#1}}

\providecommand{\citelist}[1]{#1}
\providecommand{\myeprint}[2]{E-print: \href{#1}{\texttt{#2}}}

\DeclareFontFamily{OT1}{cyr}{}
\DeclareFontShape{OT1}{cyr}{m}{n}
   {  <5> <6> <7> <8> <9> gen * wncyr
      <10> <10.95> <12> <14.4> <17.28> <20.74> <24.88> wncyr10}{}
\DeclareFontShape{OT1}{cyr}{m}{it}
    {
       <5> <6> <7> <8> <9> gen * wncyi
      <10> <10.95> <12> <14.4> <17.28> <20.74> <24.88>wncyi10
      }{}
\DeclareFontShape{OT1}{cyr}{m}{ss}
    {
       <5> <6> <7> <8> wncyss8
       <9> wncy9
      <10> <10.95> <12> <14.4> <17.28> <20.74> <24.88>wncyss10
      }{}
\DeclareFontShape{OT1}{cyr}{m}{sc}
    {
       <5> <6> <7> <8> <9> <10> <10.95> <12> <14.4> <17.28> <20.74> <24.88>wncysc10
      }{}
\DeclareFontShape{OT1}{cyr}{bx}{n}
   {
       <5> <6> <7> <8> <9> gen * wncyb
      <10> <10.95> <12> <14.4> <17.28> <20.74> <24.88>wncyb10
      }{}
\DeclareTextFontCommand{\textcyr}{\fontfamily{cyr}\selectfont}
\providecommand{\cyr}{\fontfamily{cyr}\selectfont\def\cprime{\~}}
\providecommand{\cprime}{'}

\providecommand{\Space}[3][]{\ifx#2R\ifx#1e \mathbb{C}^{#3} \else
\ifx#1p \mathbb{D}^{#3} \else
\ifx#1h \mathbb{O}^{#3} \else
\ifx#1\sigma \mathbb{A}^{#3} \else
\ensuremath{\mathbb{#2}^{#3}_{#1}{}} \fi \fi \fi \fi \else
\ensuremath{\mathbb{#2}^{#3}_{#1}{}} \fi}

\providecommand{\Space}[3][]{\ensuremath{\mathbb{#2}^{#3}_{#1}{}}}
\providecommand{\SL}[1][2]{\ensuremath{\mathrm{SL}_{#1}(\Space{R}{})}}

\providecommand{\tr}{\mathop{\mathrm{tr}}}
\providecommand{\scalar}[3][\relax]{\left\langle #2,#3 
        \right\rangle\ifx#1\relax\else_{#1}\fi}
\providecommand{\modulus}[2][\relax]{\left| #2 \right|\ifx#1\relax\else_{#1}\fi}
  \providecommand{\Zbl}[1]{Zbl\href{http://www.emis.de:80/cgi-bin/zmen/ZMATH/en/zmathf.html?first=1&maxdocs=3&type=html&an=#1&format=complete}{#1}}

\providecommand{\rmi}{\mathrm{i}}

\newtheorem{thm}{Theorem}
\newtheorem{prop}[thm]{Proposition}
\newtheorem{lem}[thm]{Lemma}
\newtheorem{cor}[thm]{Corollary}

\theoremstyle{definition}
\newtheorem{defn}[thm]{Definition}
\newtheorem{example}[thm]{Example}
\newtheorem{fig}[thm]{Figure}

\theoremstyle{remark}
\newtheorem{rem}[thm]{Remark}

\providecommand{\ccrossratio}[1]{\left\langle #1 \right\rangle}

\usepackage[msc-links,nobysame,abbrev]{amsrefs}

\BibSpecAlias{incollection}{inproceedings}
\BibSpec{article}{%
    +{}  {\PrintAuthors}                {author}
    +{.} { }                            {title}
    +{.} { }                            {part}
    +{:} { }                            {subtitle}
    +{.} { \PrintContributions}         {contribution}
    +{.} { \PrintPartials}              {partial}
    +{.} { \emph}                       {journal}
    +{,} { \textbf}                            {volume}
    +{}  { \parenthesize}               {number}
    +{:} {}                             {pages}
    +{,} { \PrintDateB}                 {date}
    +{,} { }                            {status}
    +{.} { \PrintTranslation}           {translation}
    +{.} { Reprinted in \PrintReprint}  {reprint}
    +{.} { }                            {note}
    +{.} {}                             {transition}
}

\BibSpec{partial}{%
    +{}  {}                             {part}
    +{:} { }                            {subtitle}
    +{.} { \PrintContributions}         {contribution}
    +{.} { \emph}                       {journal}
    +{,} { \textbf}                            {volume}
    +{}  { \parenthesize}               {number}
    +{:} {}                             {pages}
    +{,} { \PrintDateB}                 {date}
}

\BibSpec{book}{%
    +{}  {\PrintPrimary}                {transition}
    +{.} { \emph}                       {title}
    +{.} { }                            {part}
    +{:} { \emph}                       {subtitle}
    +{.} { }                            {series}
    +{,} { \voltext}                    {volume}
    +{.} { Edited by \PrintNameList}    {editor}
    +{.} { Translated by \PrintNameList}{translator}
    +{.} { \PrintContributions}         {contribution}
    +{.} { }                            {publisher}
    +{.} { }                            {organization}
    +{,} { }                            {address}
    +{,} { \PrintEdition}               {edition}
    +{,} { \PrintDateB}                 {date}
    +{.} { }                            {note}
    +{.} {}                             {transition}
    +{.} { \PrintTranslation}           {translation}
    +{.} { Reprinted in \PrintReprint}  {reprint}
    +{.} {}                             {transition}
}

\BibSpec{collection.article}{%
    +{}  {\PrintAuthors}                {author}
    +{.} { }                            {title}
    +{.} { }                            {part}
    +{:} { }                            {subtitle}
    +{.} { \PrintContributions}         {contribution}
    +{.} { \PrintConference}            {conference}
    +{.} { \PrintBook}                  {book}
    +{.} { In \PrintEditorsA}              {editor}
    +{.} { \emph}                         {booktitle}
    +{,} { pages~}                      {pages}
    +{,} { }                            {publisher}
    +{,} { }                            {organization}
    +{,} { }                            {address}
    +{,} { \PrintDateB}                 {date}
    +{.} { \PrintTranslation}           {translation}
    +{.} { Reprinted in \PrintReprint}  {reprint}
    +{.} { }                            {note}
    +{.} {}                             {transition}
}

 \BibSpec{inproceedings}{%
     +{}  {\PrintAuthors}                {author}
     +{.} { }                            {title}
     +{.} { }                            {part}
     +{:} { }                            {subtitle}
     +{.} { \PrintContributions}         {contribution}
     +{.} { \PrintConference}            {conference}
     +{.} { \PrintBook}                  {book}
     +{.} { In \PrintEditorsA}              {editor}
     +{.} {  \emph}                         {booktitle}
     +{,} { pages~}                      {pages}
     +{,} { }                            {publisher}
     +{,} { }                            {organization}
     +{,} { }                            {address}
     +{,} { \PrintDateB}                 {date}
     +{.} { \PrintTranslation}           {translation}
     +{.} { Reprinted in \PrintReprint}  {reprint}
     +{.} { }                            {note}
     +{.} {}                             {transition}
 }

\BibSpec{conference}{%
    +{}  {}                        {title}
    +{}  {\PrintConferenceDetails} {transition}
}

\BibSpec{innerbook}{%
    +{.} { \emph}                       {title}
    +{.} { }                            {part}
    +{:} { \emph}                       {subtitle}
    +{.} { }                            {series}
    +{,} { \voltext}                    {volume}
    +{.} { Edited by \PrintNameList}    {editor}
    +{.} { Translated by \PrintNameList}{translator}
    +{.} { \PrintContributions}         {contribution}
    +{.} { }                            {publisher}
    +{.} { }                            {organization}
    +{,} { }                            {address}
    +{,} { \PrintEdition}               {edition}
    +{,} { \PrintDateB}                 {date}
    +{.} { }                            {note}
    +{.} {}                             {transition}
}

\BibSpec{report}{%
    +{}  {\PrintPrimary}                {transition}
    +{.} { \emph}                       {title}
    +{.} { }                            {part}
    +{:} { \emph}                       {subtitle}
    +{.} { \PrintContributions}         {contribution}
    +{.} { Technical Report }           {number}
    +{,} { }                            {series}
    +{.} { }                            {organization}
    +{,} { }                            {address}
    +{,} { \PrintDateB}                 {date}
    +{.} { \PrintTranslation}           {translation}
    +{.} { Reprinted in \PrintReprint}  {reprint}
    +{.} { }                            {note}
    +{.} {}                             {transition}
}

\BibSpec{thesis}{%
    +{}  {\PrintAuthors}                {author}
    +{,} { \emph}                       {title}
    +{:} { \emph}                       {subtitle}
    +{.} { \PrintThesisType}            {type}
    +{.} { }                            {organization}
    +{,} { }                            {address}
    +{,} { \PrintDateB}                 {date}
    +{.} { \PrintTranslation}           {translation}
    +{.} { Reprinted in \PrintReprint}  {reprint}
    +{.} { }                            {note}
    +{.} {}                             {transition}
}

\makeatletter
  \def\p@enumi{\thethm}

\makeatother

\begin{document}
\title{Cycles Cross Ratio: an Invitation}

\author[Vladimir V. Kisil]%
{\href{http://www.maths.leeds.ac.uk/~kisilv/}{Vladimir V. Kisil}}
\thanks{On  leave from Odessa University.}

\address{%
School of Mathematics\\
University of Leeds\\
Leeds LS2\,9JT\\
UK
}

\email{\href{mailto:V.Kisil@leeds.ac.uk}{V.Kisil@leeds.ac.uk}}
\email{\href{mailto:kisill@maths.leeds.ac.uk}{kisilv@maths.leeds.ac.uk}}

\urladdr{\url{http://www.maths.leeds.ac.uk/~kisilv/}}

\date{\today}

\begin{abstract}
  The paper introduces \emph{cycles cross ratio}, which extends the classic cross ratio of four points to various settings: conformal  geometry, Lie spheres geometry, etc. Just like its classic counterpart cycles cross ratio is a measure of anharmonicity between spheres with respect to inversion. It also provides a M\"obius invariant distance between spheres. Many further properties of cycles cross ratio awaiting their exploration. In abstract framework the new invariant can be considered in any projective space with a bilinear pairing.
\end{abstract}
\maketitle

\section{Introduction}
\label{sec:introduction}

We introduce a new invariant in \href{https://en.wikipedia.org/wiki/Lie_sphere_geometry}{Lie sphere
  geometry}~\citelist{ \cite{Cecil08a} \cite{Benz07a}*{Ch.~3} \cite{Juhasz18a}} which is analogous to the \href{https://en.wikipedia.org/wiki/Cross-ratio}{cross-ratio} of four points in projective~\citelist{\cite{Milne11a} \cite{ShafarevichRemizov13a}*{\S~9.3}  \cite{Yaglom79}*{App.~A}} and hyperbolic geometries~\citelist{\cite{Beardon95}*{\S~4.4} \cite{Pedoe95a}*{III.5}  \cite{Schwerdtfeger79a}*{\S~I.5}}. To avoid technicalities and stay visual we will work in two dimensions, luckily all main features are already illuminated in this case. Our purpose is not to give an exhausting presentation (in fact, we are hoping it is far from being possible now), but rather to draw attention to the new invariant and its benefits to the various geometrical settings. There are several far-reaching generalisations in higher dimensions and non-Euclidean metrics, see Sect.~\ref{sec:disc-gener} for a discussion.  In abstract framework the new invariant can be considered in any projective space with a bilinear pairing.

 There is an \href{https://github.com/vvkisil/Cycles-cross-ratio-Invitation}{interactive version of this paper} implemented as a Jupyter notebook~\cite{Kisil21d} based on the \textsf{MoebInv} software package~\cite{Kisil20a}.

\section{Preliminaries: Projective space of cycles}
\label{sec:proj-spece-cycl}

We introduce our topic following its historic development. This conforms to our hands-on approach as an opposition to a start from the most general abstract construction.

Circles are the traditional subject of geometrical studies and their numerous properties were already widely presented in Euclid's \emph{Elements}. However, advances in their analytical presentations is not a so distant history.

A straightforward parametrisation of a circle equation:
\begin{equation}
  \label{eq:circles-Pedoe}
  x^2+y^2+2gx+2fy+c=0
\end{equation}
by a point \((g,f,c)\) in some subset of the three-dimensional Euclidean space \(\Space{R}{3}\) was used in~\cite{Pedoe95a}*{Ch.~II}. Abstractly, we can treat a \emph{point} \((x_0,y_0)\) of a plane as the zero radius circle with coefficients \((g,f,c)=(-x_0,-y_0,x_0^2+y_0^2)\).

It is more advantageous to use the equation
\begin{equation}
  \label{eq:cycle-generic}
  k(x^2+y^2) -2lx -2ny +m=0 \, ,
\end{equation}
which also includes \emph{straight lines} for \(k=0\). This extension comes at a price: parameters \((k,l,n,m)\) shall be treated as elements of the three dimensional projective space \(P\Space{R}{3}\) rather than the Euclidean space \(\Space{R}{4}\) since equations~\eqref{eq:cycle-generic} with \((k,l,n,m)\) and \((k_1,l_1,n_1,m_1)= (\lambda k, \lambda l, \lambda n,\lambda m)\) defines the same set of points for any \(\lambda\neq 0\). This parametrisation is known as tetracyclic/polyspheric coordinates, cf.~\citelist{\cite{Kastrup08a}*{\S~2.4.1} \cite{BergerII}*{\S~20.7}}.

The next observation is that the linear structure of \(P\Space{R}{3}\) is relevant for circles geometry. For example, the traditional concept of \emph{pencil of circles}~\cite{CoxeterGreitzer}*{\S~2.3} is nothing else but the linear span in \(P\Space{R}{3}\)~\cite{Schwerdtfeger79a}*{\S~I.1.c}. Therefore, it will be convenient to accept all points \((k,l,n,m)\in P\Space{R}{3}\) on equal ground even if they correspond to an empty set of solutions \((x,y)\) in~\eqref{eq:cycle-generic}. The latter can be thought as ``circles with imaginary radii''.

It is also appropriate to consider \((0,0,0,1)\) as a representative of the point \(  C_\infty\) at infinity, which complements \(\Space{R}{2}\) to the Riemann sphere.
Following~\cite{Yaglom79}, we call circles (with real and imaginary radii), straight lines and points (including \(  C_\infty\)) on a plane by joint name \emph{cycles}. Correspondingly, the space \(P\Space{R}{3}\) representing them---the \emph{cycles space}.

Sometimes, we need to consider cycles with an \emph{orientation}. This can be used, for example, to distinguish ``inner'' and ``outer'' tangency of cycles~\cite{FillmoreSpringer00a,Kisil20a,Kisil14b}. Furthermore, oriented cycles and their tangency relations are the starting point for the Lie sphere geometry. It is easy to encode cycles' orientations through polyspheric coordinates: cycles \((k,l,n,m)\) and \((k',l',n',m')\) have the opposite orientations if \((k,l,n,m) = (\lambda k', \lambda l', \lambda n', \lambda m')\)
for some negative \(\lambda \) and two points \((k,l,n,m)\) and \((k',l',n',m')\) are considered as representatives of the same oriented cycle if \((k,l,n,m) = (\lambda k', \lambda l', \lambda n', \lambda m')\)
for a positive \(\lambda \).

An algebraic consideration often benefits from the introduction of complex numbers. For example, we can re-write~\eqref{eq:cycle-generic} as:
\begin{align}
  \label{eq:circle-complex-matrix}
  \begin{pmatrix}
    -1&\overline{z}
  \end{pmatrix}
  \begin{pmatrix}
    \overline{L}&-m\\
    k&-{L}
  \end{pmatrix}
  \begin{pmatrix}
    z\\1
  \end{pmatrix}&=
                 k z\overline{z}-L\overline{z}-\overline{L}z+m\\
  \nonumber 
  & =   k(x^2+y^2) -2lx -2ny +m\,.
\end{align}
where \(z=x+\rmi y\) and \(L=l+\rmi n\). We call  \emph{Fillmore--Springer--Cnops construction} (or \emph{FSCc} for short) the association of the matrix \(C= \begin{pmatrix}
    \overline{L}&-m\\
    k&-{L}
  \end{pmatrix}
\) to a cycle with coefficients \((k,l,n,m)\)~\citelist{\cite{FillmoreSpringer90a} \cite{Cnops94a} \cite{Cnops02a}*{\S~4.1}}.

A reader may expect a more straightforward realisation of the quadratic form~\eqref{eq:cycle-generic}, cf.~\citelist{\cite{Schwerdtfeger79a}*{\S~I.1} \cite{Simon11a}*{Thm.~9.2.11}}:
\begin{equation}
  \label{eq:circle-complex-Schwerdtfeger}
  \begin{pmatrix}
    \overline{z} & 1
  \end{pmatrix}
  \begin{pmatrix}
    k&-{L}\\
    -\overline{L}&m
  \end{pmatrix}
  \begin{pmatrix}
    z\\1
  \end{pmatrix}=
  k z\overline{z}-L\overline{z}-\overline{L}z+m\,.
\end{equation}
However, FSCc will show its benefits in the next section. Meanwhile, we note that a point \(z\) corresponding to the zero radius cycle with the centre \(z\) is represented by the matrix
\begin{equation}
  \label{eq:zero-radius-defn}
  Z=
  \begin{pmatrix}
    \overline{z} & -z \overline{z}\\
    1 & -z
  \end{pmatrix}
  = \frac{1}{2}
  \begin{pmatrix}
    \overline{z} & \overline{z}\\
    1 & 1
  \end{pmatrix}
  \begin{pmatrix}
    1 & -z \\
    1 & -z
  \end{pmatrix} \quad \text {for} \quad  (k,l,n,m)=(1,x,y,x^2+y^2)
\end{equation}
where \(\det Z=0\). Also, the point at infinity can be represented by a zero radius cycle:
\begin{equation}
  \label{eq:infinity-cycle}
  C_\infty =
  \begin{pmatrix}
    0& -1\\0&0
  \end{pmatrix} \qquad \text{for}\qquad
  (k,l,n,m)=(0,0,0,1)\,.
\end{equation}
  More generally, \(\det   \begin{pmatrix}
    \overline{L}&-m\\
    k&-{L}
  \end{pmatrix} = -k^2 r^2
\) for \(k\neq 0\) and the cycle's radius \(r\).

\section{Fractional linear transformations and the invariant product}
\label{sec:fract-line-transf}
In the spirit of the \href{https://en.wikipedia.org/wiki/Erlangen_program}{Erlangen programme} of Felix Klein (greatly influenced by Sophus Lie), a consideration of cycle geometry is based on a group of transformations preserving this family~\cites{Kisil06a,Kisil12a}. 
Let \(M=
\begin{pmatrix}
  \alpha  &\beta\\ \gamma & \delta 
\end{pmatrix}\) be an invertible \(2\times 2\) complex matrix. Then the \emph{fractional linear transformation} (\emph{FLT} for short) of the extended complex plane \(\dot{\Space{C}{}}= \Space{C}{}\cup \{\infty\}\) is defined by:
\begin{equation}
  \label{eq:FLT-defn}
  \begin{pmatrix}
    \alpha  &\beta\\ \gamma & \delta 
  \end{pmatrix}: \ 
  z \mapsto \frac{\alpha z +\beta}{\gamma z +\delta} \, .
\end{equation}
It will be convenient to introduce the notation \(\overline{M}\) for the matrix \(
\begin{pmatrix}
  \overline{\alpha}  &\overline{\beta}\\ \overline{\gamma} & \overline{\delta} 
\end{pmatrix}\) with complex conjugated entries of \(M\). For a cycle \(C\) the matrix \(\overline{C}\) corresponds to the reflection of \(C\) in the real axis \(y=0\).  Also, due to special structure of FSCc matrix we easily check that
\begin{equation}
  \label{eq:conjugated-inverse}
  \overline{C}C={C}\overline{C}=-\det(C)I
  \qquad \text{or} \qquad
  \overline{C} \sim C^{-1} \text{ projectively if } \det(C)\neq 0.
\end{equation}

If a cycle \(C\) is composed from a non-empty set of points in \((x,y) \in \Space{R}{2}\) satisfying~\eqref{eq:cycle-generic}, then their images under  transformation~\eqref{eq:FLT-defn} form again a cycle \(C_1\) with notable link to FSCc:
\begin{lem}
  \label{le:FLT-cycle-FSCC}
  Transformation~\eqref{eq:FLT-defn} maps a cycle \(C=  \begin{pmatrix}
    \overline{L}&-m\\
    k&-{L}
  \end{pmatrix}
\) into a cycle \(C_1=  \begin{pmatrix}
    \overline{L}_1&-m_1\\
    k_1&-{L}_1
  \end{pmatrix}
  \) such that \(C_1= \overline{M} C M^{-1}\).
\end{lem}
Clearly, the map \(C \rightarrow \overline{M} C M^{-1}\) is meaningful for any cycle, including imaginary ones, thus we regard it as FLT action on the cycle space. 
For the matrix form~\eqref{eq:circle-complex-Schwerdtfeger} the above identity \(C_1= \overline{M} C M^{-1}\) needs to be replaced by the matrix congruence  \(C_1= M^* C M\)~\citelist{\cite{Schwerdtfeger79a}*{\S~II.6.e} \cite{Simon11a}*{Thm.~9.2.13}}. This difference is significant in view of the following definition.
\begin{defn}
  \label{de:cycle-product}
    For two cycles \(C\) and \(C_1\) define the \emph{cycles product} by:
  \begin{equation}
    \label{eq:cycles-product-defn}
    \scalar{C}{C_1} = -\tr(C\overline{C}_1) \, ,
  \end{equation}
  where \(\tr\) denotes the trace of a matrix.
  
  We call two cycles \(C\) and \(C_1\) \emph{orthogonal} if \(\scalar{C}{C_1}=0\).
\end{defn}
It is easy to find the explicit expression of the cycle product~\eqref{eq:cycles-product-defn}:
\begin{equation}
  \label{eq:cycle-product-explicit}
  \scalar{C}{C_1} = km_1+k_1m-2ll_1-2nn_1\,,
\end{equation}
and observe that it is linear in coefficients of the cycle \(C\) (and \(C_1\) as well).
On the other hand, it is the initial definition~\eqref{eq:cycles-product-defn}, which allows us to use the invariance of trace under matrix similarity to conclude the following.
\begin{cor}
  \label{co:cycle-product-FLT-invariant}
  The cycles product is invariant under the transformation \(C\mapsto    \overline{M} C M^{-1}\). Therefore FLT~\eqref{eq:FLT-defn} preserves orthogonality of cycles.
\end{cor}
The cycle product is a rather recent addition to the cycle geometry, see independent works \citelist{\cite{FillmoreSpringer90a} \cite{Cnops94a} \cite{Cnops02a}*{\S~4.1} \cite{Kirillov06}*{\S~4.2}}. Interestingly, expression \eqref{eq:cycles-product-defn} essentially repeats the GNS-construction in \(C^*\)-algebras~\cite{Arveson76} which is older by half of century at least.
\begin{example} \citelist{\cite{Kisil12a}*{Ch.~6} \cite{Kisil19a} \cite{Kisil14b}}
  \label{ex:orthonality-meaning}
The cycles product and cycles orthogonality encode a great amount of geometrical characteristics. For example, for cycles represented by non-empty sets of points in \(\Space{R}{2}\) we note the following:
\begin{enumerate}
\item  \label{item:quadric-flat}
  A cycle is a straight line if it is orthogonal
    \(\scalar{C}{C_\infty}=0\) to the zero radius cycle at
    infinity \(C_\infty\)~\eqref{eq:infinity-cycle}.
\item \label{it:point-zero-radius}
  A cycle \(Z\) represents a point if \(Z\) is
  self-orthogonal (\emph{isotropic}):
  \(\scalar{Z}{Z}=0\). More generally, \eqref{eq:conjugated-inverse} implies:
\begin{equation}
  \label{eq:cycles-self-product-is-det}
  \scalar{C}{C}=2\det(C)\,.
\end{equation}
\item A cycle \(C\) passes a point \(Z\) if they are orthogonal \(\scalar{C}{Z}=0\).
\item \label{it:lobachevski-line}
  A cycle \(C\) represents a line in Lobachevsky
  geometry\cite{ShafarevichRemizov13a}*{Ch.~12} if it is orthogonal
  \(\scalar{C}{C_{\Space{R}{}}}=0\)  to the real line cycle
  \(C_{\Space{R}{}}=
  \begin{pmatrix}
    \rmi &0\\0& \rmi
  \end{pmatrix}
\).
\item Two cycles are orthogonal as subsets of a plane (i.e. they have perpendicular tangents at an intersection point)  if they are
  orthogonal in the sense of Defn.~\ref{de:cycle-product}.
\item Two cycles \(C\) and \(C_1\)
  are \emph{tangent} (i.e. have a unique point in common) if
  \begin{displaymath}
    \scalar{C}{C_1}^2
    =  \scalar{C}{C}
    \scalar{C_1}{C_1}\,.
  \end{displaymath}
\item
  \label{it:inversive-distance}
  \emph{Inversive distance}~\cite{CoxeterGreitzer}*{\S~5.8} \(\theta\) of two (non-isotropic) cycles is defined by the formula:
  \begin{equation}
    \label{eq:inversive-distance}
    \theta     =  \frac{\scalar{C}{C_1}}{\sqrt{
    \scalar{C}{C} \scalar{C_1}{C_1}}}.
  \end{equation}
  In particular, the above discussed orthogonality corresponds to
  \(\theta=0\) and the tangency to \(\theta=\pm1\). For intersecting
  cycles \(\theta\) is the cosine of the intersecting
  angle. 
\item A generalisation of \emph{Steiner power} \(d(C,C_1)\) of two cycles is defined
  as, cf.~\cite{FillmoreSpringer00a}*{\S~1.1}:
  \begin{equation}
    \label{eq:steiner-power}
    d(C,C_1)=    \scalar{C}{C_1}
    + \sqrt{\scalar{C}{C}
      \scalar{C_1}{C_1}}\,,
  \end{equation}
  where both cycles \(C\) and \(C_1\) are scaled to have \(k=1\) and \(k_1=1\). Geometrically, the
  generalised Steiner power for spheres provides the square of
  \emph{tangential distance}.
\end{enumerate}
\end{example}

\begin{rem}
The cycles product is \emph{indefinite}, see~\cite{GohbergLancasterRodman05a} for an  account of the theory with some refreshing differences to the more familiar situation of inner product spaces. One illustration is the presence of self-orthogonal non-zero vectors, see Ex.~\ref{it:point-zero-radius} above. Another noteworthy observation is that the product~\eqref{eq:cycle-product-explicit} has the Lorentzian signature \((1,3)\) and \(\Space{R}{4}\) with this product is isomorphic to Minkowski space-time~\cite{Kirillov06}*{\S~4.2}. 
\end{rem}

\section{Cycles cross ratio}
\label{sec:cycles-cross-ratio}

Due to the projective nature of the cycles space (i.e. matrices \(C\) and \(\lambda C\) correspond to the same cycle) a \emph{non-zero} value of the cycle product~\eqref{eq:cycles-product-defn} is not directly meaningful. Of course, this does not affect the cycles orthogonality. A partial remedy in other cases is possible through various \emph{normalisations}~\cite{Kisil12a}*{\S~5.2}. Usually they are specified by either  of the following  conditions
\begin{enumerate}
\item \(k=1\), which is convenient for metric properties of cycles. It brings us back to the initial equation~\eqref{eq:circles-Pedoe} and is not possible for straight lines; or
\item \(\scalar{C}{C}=\pm 1\) which was suggested in \cite{Kirillov06}*{\S~4.2} and is useful, say, for tangency but is not possible for points. 
\end{enumerate}

Recall, that the projective ambiguity is elegantly balanced in the \emph{cross ratio} of four points~\citelist{\cite{Beardon95}*{\S~4.4} \cite{Pedoe95a}*{III.5}  \cite{Schwerdtfeger79a}*{\S~I.5}}:
\begin{equation}
  \label{eq:cross-ratio-points}
  (z_1,z_2;z_3,z_4) = \frac{z_1-z_3}{z_1-z_4} :  \frac{z_2-z_3}{z_2-z_4}.
\end{equation}
We use this classical pattern in the following definition.
\begin{defn}
  A \emph{cycles cross ratio}  of four cycles \(C_1\), \(C_2\), \(C_3\) and  \(C_4\) is:
  \begin{equation}
    \label{eq:cycles-cross-ratio-defn}
    \ccrossratio{C_1, C_2; C_3, C_4}= \frac{\scalar{C_1}{C_3}}{\scalar{C_1}{C_4}}:\frac{\scalar{C_2}{C_3}}{\scalar{C_2}{C_4}} 
  \end{equation}
  assuming \(\scalar{C_1}{C_4} \scalar{C_2}{C_3} \neq 0\). If \(\scalar{C_1}{C_4} \scalar{C_2}{C_3} = 0\) but \(\scalar{C_1}{C_3} \scalar{C_2}{C_4} \neq 0\) we put \(\ccrossratio{C_1, C_2; C_3, C_4} = \infty\). The cycles cross ratio is generally undefined in the remaining case of an indeterminacy \(\frac{0}{0}\). 
\end{defn}
Note that  some additional geometrical reasons may help to resolve the last situation, see the consideration of orthogonality/tangency with zero radius cycle in Ex.~\ref{ex:tangency-zero-radius}.

As an initial justification of the definition we list the following properties.
\begin{prop}
  \begin{enumerate}
  \item The cycles cross ratio is a well-defined FLT-invariant of quadruples of cycles.
  \item The cycles cross ratio of four zero radius cycles is the squared modulus of the cross ratio for the respective points:
    \begin{equation}
      \label{eq:cross-ratio-points-eq}
      \ccrossratio{Z_1, Z_2; Z_3, Z_4}=  \modulus{(z_1,z_2;z_3,z_4)}^2.
    \end{equation}
  \item There is the \emph{cancellation formula}:
    \begin{equation}
      \label{eq:cross-ratio-cancellation}
          \ccrossratio{C_1, C; C_3, C_4}    \ccrossratio{C, C_2; C_3, C_4}    = \ccrossratio{C_1, C_2; C_3, C_4} \, .
    \end{equation}
  \end{enumerate}
\end{prop}
\begin{proof}
  The first statement follows from Cor.~\ref{co:cycle-product-FLT-invariant} and the construction of the cycles cross ratio. To show the second statement we derive from~\eqref{eq:zero-radius-defn} and~\eqref{eq:cycle-product-explicit} that:
  \begin{displaymath}
    \scalar{Z_i}{Z_j} = \modulus{z_i - z_j}^2,
  \end{displaymath}
  if we use representations of zero radius cycles \(Z_i\) and \(Z_j\) by coefficients with \(k_i=k_j=1\). This implies~\eqref{eq:cross-ratio-points-eq}. A demonstration of~\eqref{eq:cross-ratio-cancellation} is straightforward.
\end{proof}

To demonstrate that there is more than just a formal similarity between the two we briefly list some applications. First, we rephrase Ex.~\ref{it:inversive-distance} in new terms.
\begin{example}
  \label{ex:inversive-dist-cross-ratio}
  The \emph{capacitance} \(\mathrm{cap}(C,C_1)\) of two cycles~\citelist{\cite{HungerbuhlerKusejko18a}*{\S~5.1} \cite{HungerbuhlerVilliger21a}*{Defn.~3}} coincides with the following cycles cross product:
    \begin{displaymath}
      \mathrm{cap}(C,C_1) = \ccrossratio{C,C_1;C_1,C} = \theta^2,
    \end{displaymath}
  where \(\theta\) is the inversive distance~\eqref{eq:inversive-distance}.
 Thereafter, FLT-invariance of the cycles cross ratio implies that the intersection angle of cycles is FLT-invariant. 
    In particular cycles are
  \begin{enumerate}
  \item \label{item:orth-if-ccrossr}
    orthogonal if \(\ccrossratio{C,C_1;C_1,C}=0\);
  \item \label{it:tangency-from-cross-ratio}
    tangent if \(\ccrossratio{C,C_1;C_1,C}=\pm 1\); and
  \item disjoint if  \(\modulus{\ccrossratio{C,C_1;C_1,C}}> 1\).
  \end{enumerate}
  Relation~\ref{item:orth-if-ccrossr} is merely a consequence of the first-order orthogonality relation \(\scalar{C}{C_1}=0\), which is fundamental to conformal and incidence geometries, cf. Ex.~\ref{ex:orthonality-meaning}(i)--(v). Meanwhile,  the tangency condition~(ii) is genuinely quadratic and shall be equally significant in Lie spheres geometry, the Steiner's porism~\cites{HungerbuhlerKusejko18a,HungerbuhlerVilliger21a}, and other questions formulated purely in terms of cycles'  tangency. 
\end{example}
\begin{example}
  \label{ex:tangency-zero-radius}
  If a non-zero radius cycle passes a zero radius cycle (point) their cross ratio has an indeterminacy $\frac{0}{0}$. Geometrically their relation can be seen in either ways: as orthogonality or tangency. Therefore, the indeterminacy of the cycle cross ratio can be geometrically resolved differently either to $0$ (indicates orthogonality) or $1$ (corresponds to tangency). More specifically (see the supporting symbolic computations in the notebook \cite{Kisil21d}):
  \begin{itemize}
  \item \emph{Orthogonality}. Consider a cycle \(Z_t\) with a fixed centre and a variable squared radius $t$. Take a generic cycle  \(C\) orthogonal  to \(Z_t\). To resolve an indeterminacy $\frac{0}{0}$ we use l'Hospital's rule at the point  $t=0$, which corresponds to \(Z_t\) becoming a zero radius cycle. This produces \(\ccrossratio{C,Z_t;Z_t,C}|_{t=0}=0\).
  \item \emph{Tangency}. For a zero radius cycle \(Z\) and passing it cycle \(C\), consider a generic cycle \(C_t=(1-t)\cdot Z+ t\cdot C\), \(t\in[0,1]\) from the pencil (linearly) spanned by \(Z\) and \(C\). Since \(C_t\) touches \(C\) Ex.~\ref{it:tangency-from-cross-ratio} implies that $\ccrossratio{C_t, C; C, C_t}=1$ for \(t > 0\) and \(C_t\) coincides with \(Z\) for \(t=0\). Thus, we can extend the value \(1\) to \([Z, C; C, Z]\) by continuity.
  \end{itemize}
  The last technique makes the cycles cross ratio meaningful for Lie spheres geometry, which extends FLT by non-point Lie transformations, when a non-zero radius cycle is sent to a zero radius one.
\end{example}
\begin{example}
  The Steiner power~\eqref{eq:steiner-power} can be written as:
  \begin{equation}
    \label{eq:steiner-cross-ratio}
    d(C,C_1)=  \ccrossratio{C,C_{\Space{R}{}};C_1,C_{\Space{R}{}}} +
    \sqrt{\ccrossratio{C,C_{\Space{R}{}};C,C_{\Space{R}{}}}} \cdot
    \sqrt{\ccrossratio{C_1,C_{\Space{R}{}};C_1,C_{\Space{R}{}}}}\, ,
  \end{equation}
  where \(C_{\Space{R}{}}\) is the real line and cycles \(C\) and \(C_1\) do not need to be normalised in any particular way. Thereafter, the Steiner power is an invariant of two cycles \(C\) and \(C_1\)  under M\"obius transformations,  since they fix the real line \(C_{\Space{R}{}}\). The M\"obius invariance is not so obvious from  expression~\eqref{eq:steiner-power}.
\end{example}

The next two applications will generalise the main features of the traditional cross ratio. Recall the other name of the cross ratio---the anharmonic ratio. The origin of the latter is as follows. Two points \(z_1\) and \(z_2\) on a line define a one-dimensional sphere with the centre \(O=\frac{1}{2}(z_1+z_2)\), which can be taken as the origin. Two points \(c_1\) and \(c_2\) are called \emph{harmonically conjugated} (with respect to \(z_1\) and \(z_2\)) if:  
\begin{displaymath}
  c_1\cdot c_2 = - z_1\cdot z_2, \qquad \text{cf. } \quad
  \includegraphics[scale=.8]{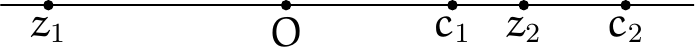}
\end{displaymath}
It is easy to check that in this case
\begin{equation}
  \label{eq:harmonic-points-cross-ratio}
  (c_1,c_2;z_1,z_2)=-1\,.
\end{equation}
Thus, the cross ratio can be viewed as a measure how far four points are from harmonic conjugation, i.e. a measure of anharmonicity of a quadruple.
\begin{figure}
  \centering
  \includegraphics[scale=.8]{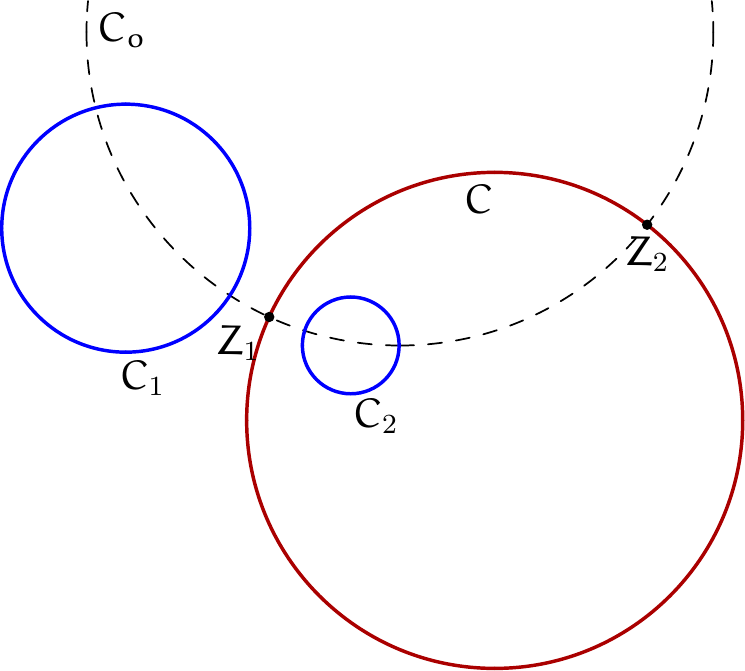}
  \caption{Cycles cross ratio for two  conjugated cycles}
  \label{fig:reflected-cycles}
\end{figure}
To make a similar interpretation of the cycles cross ratio recall that for FSCc matrices of a cycle \(C_1\) and its \emph{reflection} \(C_2\) in a cycle \(C\) we have: \(C_2=C\overline{C}_1C\), cf.~\cite{Kisil12a}*{\S~6.5}. That is, the reflection in a cycle \(C\) is the composition  of FLT transform with FSCc matrix \(C\) and complex conjugation of matrix entries. It is easy to obtain the following:
\begin{prop}
  \label{pr:inversive-distance-orthogonal}
  If a cycle \(C_1\) is a reflection of \(C_2\) in a cycle \(C\) then:
  \begin{displaymath}
    \ccrossratio{C_1, C; C, C_1} =   \ccrossratio{C_2, C; C, C_2}.
  \end{displaymath}
   More generally, the reflection in a cycle preserves the inversive distance, cf. Ex.~\ref{ex:inversive-dist-cross-ratio}. 
\end{prop}

The above condition is necessary, we describe a sufficient one as a \emph{figure} in the sense of \citelist{\cite{Kisil19a} \cite{Kisil14b} \cite{Kisil20a}}. In short, a figure is an ensemble of cycles interrelated by cycles' relations. For the purpose of this paper an FLT-invariant relation ``to be orthogonal'' between two cycles is enough. Software implementations of these figures can be found in~\cite{Kisil21d}.
\begin{fig}
\begin{enumerate}
\item For two given cycles \(C\) and \(C_1\) construct the reflection \(C_2=C\overline{C}_1C\) of \(C_1\) in \(C\).
\item Take any cycle \(C_o\) orthogonal to   \(C\) and \(C_1\), see Ill.~\ref{fig:reflected-cycles}. All such cycles make a \emph{pencil}---one dimensional subspace of the projective space of cycles. That is because there are only two linear equations for orthogonality~\eqref{eq:cycle-product-explicit} to determine four projective coordinates \((k_o,l_o,n_o,m_o)\). By Prop.~\ref{pr:inversive-distance-orthogonal} \(C_o\) is also orthogonal to \(C_2\).
\item Define cycles by orthogonality to \(C\), \(C_o\) and itself (zero radius condition), that is the intersection points of  \(C\) and \(C_o\). Since self-orthogonality is a quadratic condition there are two solutions: \(Z_1\) and \(Z_2\).
\item The harmonic conjugation of \(C_1\) and \(C_2\) (their reflection in \(C\)) implies, cf.~\eqref{eq:harmonic-points-cross-ratio}:
  \begin{equation}
    \label{eq:harmonic-cycless-cross-ratio}
    \ccrossratio{C_1,C_2;Z_1,Z_2}=1\,.
  \end{equation}
  This is demonstrated by symbolic computation in~\cite{Kisil21d}. 
\end{enumerate}
\end{fig}

\begin{figure}
  \centering
  \includegraphics[scale=.8]{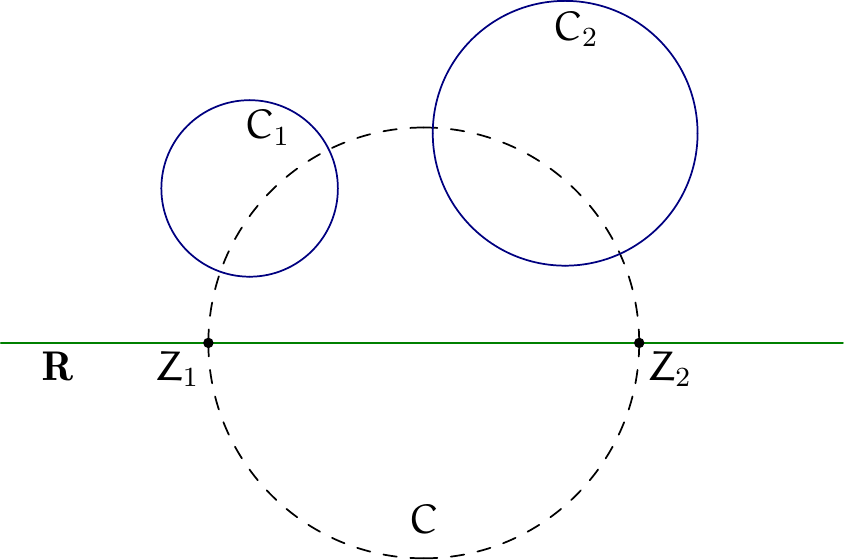}
  \caption{Construction for M\"obius invariant distance between two cycles}
  \label{fig:constr-mobi-invar}
\end{figure}
As the final illustration, we introduce M\"obius-invariant distance between cycles. Recall that FLT with an \(\SL\) matrix fixes the real line and it is called a \emph{M\"obius transformation}~\citelist{\cite{Kisil06a} \cite{Kisil12a}*{Ch.~1} \cite{Simon11a}*{\S~9.3}}. The corresponding figure is as follows, see Ill.~\ref{fig:constr-mobi-invar}:
\begin{fig}
  \label{fig:distance}
\begin{enumerate}
\item Let two distinct cycles \(C_1\) and \(C_2\)  be given and they are different from the real line \(C_{\Space{R}{}}\).
 \item Define a cycle \(C\) to be orthogonal to \(C_1\), \(C_2\), and \(C_{\Space{R}{}}\). It is specified by three linear equations for homogeneous coordinates \((k, l, n, m)\). In generic position a solution is unique, however it can be an imaginary cycle (with a negative square of the radius). 
 \item Define cycles by orthogonality to \(C\), \(C_{\Space{R}{}}\) and itself, that is the intersection points of \(C\) and \(C_{\Space{R}{}}\). In general position there are two solutions which we denote by \(Z_1\) and \(Z_2\). For an imaginary cycle \(C\) first coordinates of \(Z_1\) and \(Z_2\) are conjugated complex numbers.
 \item Since the entire construction is completely determined by the given cycles  \(C_1\) and \(C_2\) we define the \emph{distance between two cycles} by:
\begin{equation}
  \label{eq:distance-cross-ratio}
  d(C_1,C_2)= \textstyle \frac{1}{2}\log \ccrossratio{C_1,C_2; Z_1,Z_2}.
\end{equation}
\end{enumerate}
\end{fig}

From M\"obius invariance of the real line and cycles cross ratio our construction implies the following:
\begin{prop}
  \label{pr:inv-distance}
  \begin{enumerate}
  \item   Distance~\eqref{eq:distance-cross-ratio} is M\"obius invariant;
  \item \label{item:distance-equality}
    For zero radius cycles \(C_1\), \(C_2\) formula~\eqref{eq:distance-cross-ratio} coincides with the Lobachevsky metric on the upper-half plane. 
  \item \label{item:additivity}
    For any cycle \(C_3\) orthogonal to \(C\) the (signed) distance is additive: \(d(C_1,C_3)=d(C_1,C_2)+d(C_2,C_3)\).
  \item \label{item:vertical-geodesic}
    If centres of \(C_1\) and \(C_2\) are on the imaginary axis (therefore \(Z_1\) and \(Z_2\) are zero and infinity) then \(  d(C_1,C_2)=\log(m_1/k_1)-\log(k_2/m_2)\).
  \end{enumerate}
\end{prop}
Obviously, Prop.~\ref{item:distance-equality} is the consequence of~\eqref{eq:cross-ratio-points-eq} and~\ref{item:additivity} follows from the cancellation rule~\eqref{eq:cross-ratio-cancellation}. On the other hand, the expression in~\ref{item:vertical-geodesic} can be obtained by a direct computation, see~\cite{Kisil21d}. Note, that in this case \(m_i/k_i\) is the square \emph{tangential distance} (also known as the generalised Steiner power~\eqref{eq:steiner-power}, \eqref{eq:steiner-cross-ratio}) from \(Z_1\) (the origin) to the cycle \(C_i\). Thus, for a zero radius \(C_i\) it coincides with the usual distance between centres of \(Z_1\) and \(C_i\) and Prop.~\ref{item:vertical-geodesic} recreates the classical result~\cite{Beardon95}*{(7.2.6)}.

Comparing constructions on Ill.~\ref{fig:reflected-cycles}--\ref{fig:constr-mobi-invar} and formulae~\eqref{eq:harmonic-cycless-cross-ratio}--\eqref{eq:distance-cross-ratio} we can say that the invariant distance measures how far two cycles are from being reflections of each other in the real line.

\section{Discussion and generalisations}
\label{sec:disc-gener}

We presented some evidence that the cycles cross ratio extends to Lie spheres geometry the concept of the cross ratio of four points. It is natural to expect that a majority of the classic theory~\citelist{\cite{Milne11a} \cite{ShafarevichRemizov13a}*{\S~9.3}  \cite{Yaglom79}*{App.~A} \cite{Beardon95}*{\S~4.4} \cite{Pedoe95a}*{III.5}  \cite{Schwerdtfeger79a}*{\S~I.5}} admits similar adaptation as well. However, we can expect even more than that.

One can lay down a general framework for the introduced invariant~\eqref{eq:cycles-cross-ratio-defn} in generic projective spaces as follows\footnote{I am grateful to an anonymous referee for pointing this out in his/her report.}. Let \(V\) be a vector space over an arbitrary field \(\Space{F}{}\) with a bilinear pairing \(\scalar{\cdot}{\cdot}: V \times  V \rightarrow \Space{F}{}\). Upon choosing any four points of the projective space on \(V\) one may select
arbitrary non-zero vectors, say \(c_1\), \(c_ 2\), \(c_3\), \(c_4\), representing these points. Then the quotient
\begin{equation}
  \label{eq:generic-cross-ratio}
  \frac{\scalar{c_1}{c_3}}{\scalar{c_1}{c_4}} : \frac{\scalar{c_2}{c_3}}{\scalar{c_2}{c_4}}  
\end{equation}
will in general be an element of \(\Space{F}{}\). Whenever this scalar is well defined
it obviously will be an invariant under the natural action of the general
orthogonal group \(\mathrm{GO}(V, \scalar{\cdot}{\cdot})\) on the point set of the underlying projective space. More generally, a sort of projective space can be defined on a module \(V\) for a weaker structure than a field, say, algebras of dual and double numbers~\cite{Mustafa17a}. In such cases numerous divisors of zero prompt a projective treatment~\citelist{ \cite{Brewer12a} \cite{Kisil12a}*{\S~4.5} } of the new cross ratio~\eqref{eq:generic-cross-ratio}. Another perspective direction to research are discrete M\"obius  geometries~\cites{HungerbuhlerKusejko18a,HungerbuhlerVilliger21a}.

The geometric applications of the new invariant~\eqref{eq:cycles-cross-ratio-defn}, \eqref{eq:generic-cross-ratio} are expected much beyond the currently presented situation of circles on a plane. Indeed, FSCc and cycles product based on Clifford algebras works in spaces of higher dimensions and with non-degenerate metrics of arbitrary signatures~\citelist{\cite{FillmoreSpringer90a} \cite{Cnops94a} \cite{Cnops02a}*{\S~4.1} \cite{Kisil12a}}. In a straightforward fashion cycles cross ratio~\eqref{eq:cycles-cross-ratio-defn} remains a  geometric FLT-invariant in higher dimensions as well.

A more challenging situation occurs if we have a degenerate metric and cycles are represented by parabolas~\cite{Yaglom79} and respective FLT are based on dual numbers~\cite{BoltFerdinandsKavlie09a}. Besides theoretical interest such spaces are meaningful physical models~\cite{Yaglom79,Gromov12a,Gromov10a,Kastrup08a}. The differential geometry loses its ground in the degenerate case and non-commutative/non-local effects appear~\cite{Kisil12a}*{\S~7.2}. The presence of zero divisors among dual numbers prompts a projective approach~\citelist{\cite{Brewer12a} \cite{Kisil12a}*{\S~4.5}} to the cross ratio of four dual numbers, which replaces~\eqref{eq:cross-ratio-points}.  Variational methods do not produce FLT-invariant family of geodesics in the degenerate case, instead geometry of cycles needs to be employed~\cite{Kisil08a}. Our construction from Figure~\ref{fig:distance} shall be usable to define M\"obius invariant distance between parabolic cycles as well, cf.~\cite{Kisil12a}*{\S~9.5}.

Last but not least, if we restrict the group of transformations from FLT to M\"obius maps (or any other subgroup of FLT which fixes a particular cycle) we will get a larger set of invariants. In particular, the FSCc matrix~\eqref{eq:circle-complex-matrix} and the corresponding cycles product~\eqref{eq:cycles-product-defn} can use different number systems (complex, dual or double) independently from the geometry of the plane~\citelist{ \cite{Kisil06a} \cite{Kisil12a}*{\S~5.3} \cite{Kisil15a}}. Therefore, there will be three different cycles cross ratios~\eqref{eq:cycles-cross-ratio-defn} for each geometry of circles, parabolas and hyperbolas. This echoes the existence of nine Cayley--Klein geometries~\citelist{ \cite{Yaglom79}*{App.~A} \cite{Gromov90a} \cite{Pimenov65a} }. 

Various aspects of the cycles cross ratio appear to be a wide and fruitful field for further research.

\section*{Acknowledgments}
\label{sec:acknowledgments}
I am grateful to an anonymous referee for a detailed report with useful suggestions.


\small

\begin{bibdiv}
\begin{biblist}

\bib{Arveson76}{book}{
      author={Arveson, William},
       title={An invitation to {C}*-algebras},
   publisher={Springer-Verlag},
        date={{\noopsort{}}1976},
}

\bib{Beardon95}{book}{
      author={Beardon, Alan~F.},
       title={The geometry of discrete groups},
      series={Graduate Texts in Mathematics},
   publisher={Springer-Verlag},
     address={New York},
        date={1995},
      volume={91},
        ISBN={0-387-90788-2},
        note={Corrected reprint of the 1983 original},
      review={\MR{MR1393195 (97d:22011)}},
}

\bib{Benz07a}{book}{
      author={Benz, Walter},
       title={Classical geometries in modern contexts. {Geometry} of real inner
  product spaces},
     edition={Second edition},
   publisher={Birkh\"auser Verlag},
     address={Basel},
        date={2007},
        ISBN={978-3-7643-8540-8},
         url={http://dx.doi.org/10.1007/978-3-0348-0420-2},
      review={\MR{MR2370626 (2008i:51001)}},
}

\bib{BergerII}{book}{
      author={Berger, Marcel},
       title={Geometry. {II}},
      series={Universitext},
   publisher={Springer-Verlag},
     address={Berlin},
        date={1987},
        ISBN={3-540-17015-4},
        note={Translated from the French by M. Cole and S. Levy},
      review={\MR{882916 (88a:51001b)}},
}

\bib{BoltFerdinandsKavlie09a}{article}{
      author={Bolt, Michael},
      author={Ferdinands, Timothy},
      author={Kavlie, Landon},
       title={The most general planar transformations that map parabolas into
  parabolas},
        date={2009},
        ISSN={1944-4176},
     journal={Involve},
      volume={2},
      number={1},
       pages={79\ndash 88},
         url={http://dx.doi.org/10.2140/involve.2009.2.79},
      review={\MR{2501346 (2010e:51003)}},
}

\bib{Brewer12a}{article}{
      author={Brewer, Sky},
       title={Projective cross-ratio on hypercomplex numbers},
        date={2013March},
     journal={Adv. Appl. Clifford Algebras},
      volume={23},
      number={1},
       pages={1\ndash 14},
        note={\doi{10.1007/s00006-012-0335-7}. \arXiv{1203.2554}},
}

\bib{Cecil08a}{book}{
      author={Cecil, Thomas~E.},
       title={{Lie} sphere geometry: {With} applications to submanifolds},
     edition={Second},
      series={Universitext},
   publisher={Springer, New York},
        date={2008},
        ISBN={978-0-387-74655-5},
      review={\MR{2361414 (2008h:53091)}},
}

\bib{Cnops94a}{thesis}{
      author={Cnops, Jan},
       title={{Hurwitz} pairs and applications of {M\"obius} transformations},
        type={{Habilitation} Dissertation},
     address={Universiteit Gent},
        date={1994},
        note={See also~\cite{Cnops02a}},
}

\bib{Cnops02a}{book}{
      author={Cnops, Jan},
       title={An introduction to {Dirac} operators on manifolds},
      series={Progress in Mathematical Physics},
   publisher={Birkh\"auser Boston Inc.},
     address={Boston, MA},
        date={2002},
      volume={24},
        ISBN={0-8176-4298-6},
      review={\MR{1 917 405}},
}

\bib{CoxeterGreitzer}{book}{
      author={Coxeter, H.S.M.},
      author={Greitzer, S.L.},
       title={Geometry revisited},
   publisher={Random House},
     address={New York},
        date={1967},
        note={\Zbl{0166.16402}},
}

\bib{FillmoreSpringer90a}{article}{
      author={Fillmore, Jay~P.},
      author={Springer, A.},
       title={M\"obius groups over general fields using {Clifford} algebras
  associated with spheres},
        date={1990},
        ISSN={0020-7748},
     journal={Internat. J. Theoret. Phys.},
      volume={29},
      number={3},
       pages={225\ndash 246},
         url={http://dx.doi.org/10.1007/BF00673627},
      review={\MR{1049005 (92a:22016)}},
}

\bib{FillmoreSpringer00a}{article}{
      author={Fillmore, Jay~P.},
      author={Springer, Arthur},
       title={Determining circles and spheres satisfying conditions which
  generalize tangency},
        date={2000},
        note={preprint,
  \url{http://www.math.ucsd.edu/~fillmore/papers/2000LGalgorithm.pdf}},
}

\bib{GohbergLancasterRodman05a}{book}{
      author={Gohberg, Israel},
      author={Lancaster, Peter},
      author={Rodman, Leiba},
       title={Indefinite linear algebra and applications},
   publisher={Birkh\"auser Verlag},
     address={Basel},
        date={2005},
        ISBN={978-3-7643-7349-8; 3-7643-7349-0},
      review={\MR{2186302 (2006j:15001)}},
}

\bib{Gromov90a}{book}{
      author={Gromov, N.~A.},
       title={{\cyr Kontraktsii i analiticheskie prodolzheniya klassicheskikh
  grupp. {Edinyi} podkhod}. ({Russian}) [{Contractions} and analytic extensions
  of classical groups. {Unified} approach]},
   publisher={Akad. Nauk SSSR Ural. Otdel. Komi Nauchn. Tsentr},
     address={Syktyvkar},
        date={1990},
      review={\MR{MR1092760 (91m:81078)}},
}

\bib{Gromov10a}{article}{
      author={Gromov, N.~A.},
       title={Possible quantum kinematics. {II}. {Nonminimal} case},
        date={2010},
        ISSN={0022-2488},
     journal={J. Math. Phys.},
      volume={51},
      number={8},
       pages={083515, 12},
         url={http://dx.doi.org/10.1063/1.3460841},
      review={\MR{2683557 (2011k:81161)}},
}

\bib{Gromov12a}{book}{
      author={{Gromov}, N.~A.},
       title={{\cyr Kontraktsii Klassicheskikh i Kvantovykh Grupp.}
  [{C}ontractions of classic and quanrum groups]},
    language={Russian},
   publisher={Moskva: Fizmatlit},
        date={2012},
        ISBN={978-5-9221-1398-4/hbk; 978-5-7691-2325-2/hbk},
}

\bib{HungerbuhlerKusejko18a}{article}{
      author={Hungerb\"{u}hler, Norbert},
      author={Kusejko, Katharina},
       title={{S}teiner's porism in finite {M}iquelian {M}\"{o}bius planes},
        date={2018},
        ISSN={1615-715X},
     journal={Adv. Geom.},
      volume={18},
      number={1},
       pages={55\ndash 68},
         url={https://doi.org/10.1515/advgeom-2017-0027},
      review={\MR{3750254}},
}

\bib{HungerbuhlerVilliger21a}{article}{
      author={Hungerb\"{u}hler, Norbert},
      author={Villiger, Gideon},
       title={Exotic {S}teiner chains in {M}iquelian {M}\"{o}bius planes of odd
  order},
        date={2021},
        ISSN={1615-715X},
     journal={Adv. Geom.},
      volume={21},
      number={2},
       pages={207\ndash 220},
         url={https://doi.org/10.1515/advgeom-2020-0035},
      review={\MR{4243958}},
}

\bib{Juhasz18a}{article}{
      author={Juh\'{a}sz, M\'{a}t\'{e}~L.},
       title={A universal linear algebraic model for conformal geometries},
        date={2018},
        ISSN={0047-2468},
     journal={J. Geom.},
      volume={109},
      number={3},
       pages={Art. 48, 18},
         url={https://doi.org/10.1007/s00022-018-0451-1},
        note={\arXiv{1603.06863}},
      review={\MR{3871573}},
}

\bib{Kastrup08a}{article}{
      author={Kastrup, H.~A.},
       title={On the advancements of conformal transformations and their
  associated symmetries in geometry and theoretical physics},
        date={2008},
        ISSN={1521-3889},
     journal={Annalen der Physik},
      volume={17},
      number={9--10},
       pages={631\ndash 690},
         url={http://dx.doi.org/10.1002/andp.200810324},
        note={\arXiv{0808.2730}},
}

\bib{Kirillov06}{book}{
      author={Kirillov, A.~A.},
       title={A tale of two fractals},
   publisher={Springer, New York},
        date={2013},
        ISBN={978-0-8176-8381-8; 978-0-8176-8382-5},
         url={http://dx.doi.org/10.1007/978-0-8176-8382-5},
        note={Draft:
  \url{http://www.math.upenn.edu/~kirillov/MATH480-F07/tf.pdf}},
      review={\MR{3060066}},
}

\bib{Kisil08a}{article}{
      author={Kisil, Anastasia~V.},
       title={Isometric action of {${\rm SL}_2(\mathbb{R})$} on homogeneous
  spaces},
        date={2010},
     journal={Adv. App. Clifford Algebras},
      volume={20},
      number={2},
       pages={299\ndash 312},
        note={\arXiv{0810.0368}. \MR{2012b:32019}},
}

\bib{Kisil06a}{article}{
      author={Kisil, Vladimir~V.},
       title={Erlangen program at large--0: Starting with the group {${\rm
  SL}\sb 2({\bf R})$}},
        date={2007},
        ISSN={0002-9920},
     journal={Notices Amer. Math. Soc.},
      volume={54},
      number={11},
       pages={1458\ndash 1465},
        note={\arXiv{math/0607387},
  \href{http://www.ams.org/notices/200711/tx071101458p.pdf}{On-line}.
  \Zbl{1137.22006}},
      review={\MR{MR2361159}},
}

\bib{Kisil12a}{book}{
      author={Kisil, Vladimir~V.},
       title={Geometry of {M}\"obius transformations: {E}lliptic, parabolic and
  hyperbolic actions of {$\mathrm{SL}_2(\mathbf{R})$}},
   publisher={Imperial College Press},
     address={London},
        date={2012},
        note={Includes a live DVD. \Zbl{1254.30001}},
}

\bib{Kisil15a}{article}{
      author={Kisil, Vladimir~V.},
       title={{P}oincar\'e extension of {M}\"obius transformations},
        date={2017},
     journal={Complex Variables and Elliptic Equations},
      volume={62},
      number={9},
       pages={1221\ndash 1236},
        note={\arXiv{1507.02257}},
}

\bib{Kisil14b}{article}{
      author={Kisil, Vladimir~V.},
       title={An extension of {Mobius--Lie} geometry with conformal ensembles
  of cycles and its implementation in a {GiNaC} library},
        date={2018},
        ISSN={2072-9812},
     journal={Proc. Int. Geom. Cent.},
      volume={11},
      number={3},
       pages={45\ndash 67},
         url={https://doi.org/10.15673/tmgc.v11i3.1203},
        note={\arXiv{1512.02960}. Project page:
  \url{http://moebinv.sourceforge.net/}},
}

\bib{Kisil19a}{incollection}{
      author={Kisil, Vladimir~V.},
       title={{M\"obius--Lie} geometry and its extension},
        date={2019},
   booktitle={Geometry, integrability and quantization {XX}},
      editor={Mladenov, Iva\"{\i}lo~M.},
      editor={Meng, Guowu},
      editor={Yoshioka, Akira},
   publisher={Bulgar. Acad. Sci., Sofia},
       pages={13\ndash 61},
        note={\arXiv{1811.10499}},
}

\bib{Kisil20a}{article}{
      author={Kisil, Vladimir~V.},
       title={{MoebInv}: {C++} libraries for manipulations in non-{Euclidean}
  geometry},
        date={2020},
        ISSN={2352-7110},
     journal={SoftwareX},
      volume={11},
  pages={\href{http://www.sciencedirect.com/science/article/pii/S2352711019302523}{100385}},
  url={http://www.sciencedirect.com/science/article/pii/S2352711019302523},
        note={\doi{10.1016/j.softx.2019.100385}},
}

\bib{Kisil21d}{article}{
      author={Kisil, Vladimir~V.},
       title={Cycles cross ratio: a {Jupyter} notebook},
        date={2021},
     journal={GitHub},
        note={\url{https://github.com/vvkisil/Cycles-cross-ratio-Invitation}},
}

\bib{Milne11a}{book}{
      author={Milne, John~J.},
       title={An elementary treatise on cross-ratio geometry: with historical
  notes},
   publisher={Cambridge: The University Press},
        date={1911},
  note={\url{https://archive.org/details/elementarytreati00milnuoft/page/n5/mode/2up}},
}

\bib{Mustafa17a}{article}{
      author={Mustafa, Khawlah~A.},
       title={The groups of two by two matrices in double and dual numbers, and
  associated {M\"{o}bius} transformations},
        date={2018},
        ISSN={0188-7009},
     journal={Adv. Appl. Clifford Algebr.},
      volume={28},
      number={5},
       pages={Art. 92, 25},
         url={https://doi.org/10.1007/s00006-018-0910-7},
        note={\arXiv{1707.01349}},
      review={\MR{3860130}},
}

\bib{Pedoe95a}{book}{
      author={Pedoe, Dan},
       title={Circles: A mathematical view},
      series={MAA Spectrum},
   publisher={Mathematical Association of America, Washington, DC},
        date={1995},
        ISBN={0-88385-518-6},
        note={Revised reprint of the 1979 edition, With a biographical appendix
  on Karl Feuerbach by Laura Guggenbuhl},
      review={\MR{1349339 (96e:51020)}},
}

\bib{Pimenov65a}{article}{
      author={Pimenov, R.I.},
       title={Unified axiomatics of spaces with maximal movement group},
    language={Russian},
        date={1965},
     journal={Litov. Mat. Sb.},
      volume={5},
       pages={457\ndash 486},
        note={\Zbl{0139.37806}},
}

\bib{Schwerdtfeger79a}{book}{
      author={Schwerdtfeger, Hans},
       title={Geometry of complex numbers: Circle geometry, {Moebius}
  transformation, non-{Euclidean} geometry},
      series={Dover Books on Advanced Mathematics},
   publisher={Dover Publications Inc.},
     address={New York},
        date={1979},
        ISBN={0-486-63830-8},
        note={A corrected reprinting of the 1962 edition},
      review={\MR{620163 (82g:51032)}},
}

\bib{ShafarevichRemizov13a}{book}{
      author={Shafarevich, Igor~R.},
      author={Remizov, Alexey~O.},
       title={Linear algebra and geometry},
   publisher={Springer, Heidelberg},
        date={2013},
        ISBN={978-3-642-30993-9},
         url={https://doi.org/10.1007/978-3-642-30994-6},
        note={Translated from the 2009 Russian original by David Kramer and
  Lena Nekludova},
      review={\MR{2963561}},
}

\bib{Simon11a}{book}{
      author={Simon, Barry},
       title={{Szeg\H o's} theorem and its descendants. {Spectral} theory for
  {$L^2$} perturbations of orthogonal polynomials},
      series={M. B. Porter Lectures},
   publisher={Princeton University Press, Princeton, NJ},
        date={2011},
        ISBN={978-0-691-14704-8},
      review={\MR{2743058}},
}

\bib{Yaglom79}{book}{
      author={Yaglom, I.~M.},
       title={A simple non-{Euclidean} geometry and its physical basis},
      series={Heidelberg Science Library},
   publisher={Springer-Verlag},
     address={New York},
        date={1979},
        ISBN={0-387-90332-1},
        note={Translated from the Russian by Abe Shenitzer, with the editorial
  assistance of Basil Gordon},
      review={\MR{MR520230 (80c:51007)}},
}

\end{biblist}
\end{bibdiv}

\providecommand{\noopsort}[1]{} \providecommand{\printfirst}[2]{#1}
  \providecommand{\singleletter}[1]{#1} \providecommand{\switchargs}[2]{#2#1}
  \providecommand{\irm}{\textup{I}} \providecommand{\iirm}{\textup{II}}
  \providecommand{\vrm}{\textup{V}} \providecommand{\cprime}{'}
  \providecommand{\eprint}[2]{\texttt{#2}}
  \providecommand{\myeprint}[2]{\texttt{#2}}
  \providecommand{\arXiv}[1]{\myeprint{http://arXiv.org/abs/#1}{arXiv:#1}}
  \providecommand{\doi}[1]{\href{http://dx.doi.org/#1}{doi:
  #1}}\providecommand{\CPP}{\texttt{C++}}
  \providecommand{\NoWEB}{\texttt{noweb}}
  \providecommand{\MetaPost}{\texttt{Meta}\-\texttt{Post}}
  \providecommand{\GiNaC}{\textsf{GiNaC}}
  \providecommand{\pyGiNaC}{\textsf{pyGiNaC}}
  \providecommand{\Asymptote}{\texttt{Asymptote}}

\end{document}